\documentclass[11pt, letterpaper]{amsart}

\usepackage{amsmath}
\usepackage{amsfonts}
\usepackage{tikz}
\usepackage{amsthm}
\usepackage[pdftex, colorlinks=true]{hyperref}

\theoremstyle{plain}
\newtheorem*{theorem}{Theorem}

\newcommand{\Z}{\mathbb{Z}}
\newcommand{\tp}{\texorpdfstring}
\newcommand{\symgroup}[1]{Sym_{#1}}
\newcommand{\braidgroup}[1]{Br_{#1}}
\newcommand{\spacedot}{{ \ \cdot{} \ }}

\title[Differently knotted symplectic surfaces in $D^4$]{Differently knotted symplectic surfaces in $D^4$ bounded by the same transverse knot}
\author{Andrew Geng}

\begin{document}

\maketitle

\section{Introduction}

This paper is concerned with symplectic surfaces in the four-dimensional 
ball. More precisely, we consider
embedded connected surfaces $S \subset D^4$ whose boundary $\partial S = S 
\cap \partial D^4 \subset S^3$ is a transverse knot. We prove:

\begin{theorem}
There are two symplectic surfaces $S_1$ and $S_2$ which bound the same 
transverse knot, have the same topology (as abstract surfaces), and such 
that $\pi_1(D^4 \setminus S_1)$ is not isomorphic to $\pi_1(D^4 \setminus 
S_2)$.
\end{theorem}

This builds on a previous example (\cite{Auroux}, \S5) of surfaces bounding the 
same transverse link, which had different topology (one was connected, and 
the other was not). We add the same piece to both these examples to 
construct ours.

\subsection{Acknowledgments}

This work was done under the mentorship of Paul Seidel
and supported by a grant from the Massachusetts Institute of Technology.

\section{Background}

It is well-known \cite{Bennequin} that 
the closure of any braid $\beta \in \braidgroup{m}$ is naturally a transverse link. 
This has been used to construct examples of transverse links which lie in 
the same topological isotopy class, but are not isotopic as transverse 
links \cite{BirmanMenasco}. A factorization of $\beta$ is an expression
\[
\beta = \sigma_1 \cdots \sigma_k
\]
where each $\sigma_i$ is conjugate to one of the standard Artin generators 
of $\braidgroup{m}$. Every factorization of $\beta$ describes a symplectic surface 
$S \subset D^4$ whose boundary is the transverse link associated to
$\beta$ (see e.g. \cite{AurouxSmith}).
$S$ is connected if and only if the images of
$\sigma_1,\dots,\sigma_k$ in the symmetric group $\symgroup{m}$ act transitively. 
The Euler characteristic of $S$ is $m - k$;
hence the topological type of the abstract surface $S$
depends only on $\beta$, and not on the factorization.

There is a well-known method \cite{Teicher} for computing a presentation $\pi_1(D^4 
\setminus S)$, as a quotient of the fundamental group of the $m$-punctured 
disc, which is $F_m = \langle x_1, \dots, x_m \rangle$. Every word 
$\sigma_i$ in the factorization yields a relation. This is best 
represented graphically as follows. Conjugates of the standard Artin 
generators correspond bijectively to embedded paths in the disc
(up to isotopy) joining two of the punctures.

Given any such path, one fattens it into a figure-eight loop
enclosing the punctures at its endpoints. From a path $\gamma$
joining two punctures $p_i$ and $p_j$, one obtains a loop
$\gamma'c_i(\gamma')^{-1}c_j^{-1}$, where $c_i$ and $c_j$
are counterclockwise circles around the punctures and $\gamma'$
is an appropriate segment of $\gamma$.

This identifies a conjugacy class in $F_m$, any element of which
can be taken as the relation implied by this path. For consistency
or ease of computation we may choose any convenient basepoint
to determine what this relation is.

For example, let $a$, $b$, and $c$ denote the Artin generators in $\braidgroup{4}$.
The word $(ac) b (ac)^{-1}$ yields the following:
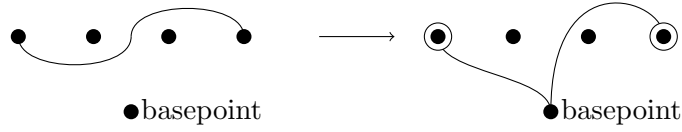
\begin{center}
\begin{tikzpicture}
	\path (2.5,-1) coordinate[label=right:basepoint] (x0);
	\fill [black] (x0) circle (0.1);

	\draw (1,0) .. controls (1,-0.5) and (2.5,-0.5) .. (2.5,0);
	\draw (2.5,0) .. controls (2.5,0.5) and (4,0.5) .. (4,0);

	\foreach \point in {(1,0), (2,0), (3,0), (4,0)}
		\fill [black] \point circle (0.1);

	\draw[->] (5,0) -- (6,0) node[right] {};
\end{tikzpicture}
\begin{tikzpicture}
	\path (2.5,-1) coordinate[label=right:basepoint] (x0);
	\fill [black] (x0) circle (0.1);

	\draw (1,0) node[draw=black, fill=white, shape=circle] {} .. controls (1,-0.5) and (2.5,-0.5) .. (x0);
	\draw (x0) .. controls (2.5,1) and (4,0.5) .. (4,0) node[draw=black, fill=white, shape=circle] {};

	\foreach \point in {(1,0), (2,0), (3,0), (4,0)}
		\fill [black] \point circle (0.1);
\end{tikzpicture}
\end{center}
\begin{align*}
	x_1 &= x_3^{-1} x_4 x_3
\end{align*}

\section{The Examples}

Let $a$, $b$, and $c$ be the Artin generators in $\braidgroup{4}$.
Appending the factor $a^3 b a^{-3}$ to the examples given in \cite{Auroux} \S5
yields the following factorizations:
\begin{align}
	& (c^{-2} b) a (c^{-2} b)^{-1} \spacedot b \spacedot (ac^3) b (ac^3)^{-1} \spacedot \\
	& \hspace{\parindent} (ac^5 b^{-1}) a (ac^5 b^{-1})^{-1} \spacedot c \spacedot c \spacedot a^3 b a^{-3} \label{one} \nonumber\\[1ex]
	& b \spacedot (ac)b(ac)^{-1} \spacedot (ac)b(ac)^{-1} \spacedot (ac)^2 b(ac)^{-2} \spacedot \\
	& \hspace{\parindent} (ac)^2 b(ac)^{-2} \spacedot (ac)^3b(ac)^{-3} \spacedot a^3 b a^{-3} \label{two} \nonumber
\end{align}

It was checked in \cite{Auroux} that these factorizations are of the same
braid. Their images in $\symgroup{4}$ are:
\begin{align}
	(13) \cdot (23) \cdot (14) \cdot (24) \cdot (34) \cdot (34) \cdot (13) \\
	(23) \cdot (14) \cdot (14) \cdot (23) \cdot (23) \cdot (14) \cdot (13)
\end{align}
The transpositions $(13)$, $(14)$, and $(23)$ are present in both and suffice
to generate $\symgroup{4}$, so $S_1$ and $S_2$ are both connected.
The product is the cyclic permutation $(1234)$, so the boundaries
$\partial S_1$ and $\partial S_2$ are also connected.

From factorization (\ref{one}), we have the following relations:
\begin{align*}
	b & & x_2 &= x_3 \\
	c & & x_3 &= x_4 \\
	(c^{-2} b) a (c^{-2} b)^{-1} & & x_1 &= x_2^{-1} x_4 x_3 x_4^{-1} x_2 = x_3
\end{align*}
Given these relations, the relations arising from the rest of the factors
simplify to identity. With, $x_1 = x_2 = x_3 = x_4$,
$\pi_1(D^4 \setminus S_1)$ must be $\Z$.

From factorization (\ref{two}), we have:
\begin{align*}
	b & & x_2 &= x_3 \\
	(ac) b (ac)^{-1} & & x_1 &= x_3^{-1} x_4 x_3 \\
	a^3 b a^{-3} & & x_3 &= x_1^{-1} x_2^{-1} x_1 x_2 x_1
\end{align*}
This last relation can be rewritten, using $x_2 = x_3$, as:
\[ x_1 x_2 x_1 = x_2 x_1 x_2 \]
As with the previous example, the remaining factors yield relations
that simplify to identity given these three.
Thus $\pi_1(D^4 \setminus S_2) = \braidgroup{3}$.

\section{A Note on Double Branched Covers of \tp{$D^4$}{the 4-ball}}

Given any connected surface $S \subset D^4$, one can form the double
cover of $D^4$ branched along $S$, which we denote by $M$. If $S$ is symplectic
and its boundary is a transverse link, $M$ is an exact symplectic
manifold with contact type boundary \cite{AurouxSmith}.

The fundamental group $\pi_1(D^4 \setminus S)$ comes with a canonical homomorphism
$\phi: \pi_1(D^4 \setminus S) \rightarrow \Z$. The meridian element $x$,
determined up to conjugacy, satisfies $\phi(x) = 1$. To compute $\pi_1(M)$,
one takes the subgroup $\pi_1(D^4 \setminus S)' = \phi^{-1}(2\Z)$,
and then divides out by $x^2$
(this is a simple application of the Seifert-van Kampen theorem).

Take the examples $S_1$, $S_2$ from the previous Section,
and let $M_1$, $M_2$ be the associated double branched covers.
It is easy to see that $\pi_1(M_1)$ is trivial.
In the second case, we have a homomorphism
$\pi_1(D^4 \setminus S_2) \cong \braidgroup{3} \rightarrow \symgroup{3}$,
which restricts to $\pi_1(D^4 \setminus S_2)' \rightarrow A_3 \cong \Z/3$.
Since $x^2$ goes to zero under this, we get an induced homomorphism
$\pi_1(M_2) \rightarrow \Z/3$. One easily checks that this is surjective.
Hence, $M_1$ and $M_2$ are two different exact symplectic fillings of the
same contact three-manifold $\partial M_1 = \partial M_2$.

For other such examples and further discussion,
see \cite{Akhmedov}, \cite{Ozbagci}, and \cite{Smith}.

\section{Closing Remarks}

Numerical evidence suggests that, for odd $s \geq 3$, using $a^s b a^{-s}$
in place of $a^3 b a^{-3}$ makes
$\pi_1(D^4 \setminus S_2) = \left< x,y \mid x^2 = y^s \right>$
while $\pi_1(D^4 \setminus S_1)$ remains as $\Z$.
These groups are distiguished from each other by the number of homomorphisms
from them to the dihedral groups.

Electronic resources (Python code)
for reproducing the numerical results can be found
online at \url{http://www-math.mit.edu/~seidel/geng/}.


\begin{thebibliography}{}

\bibitem{Akhmedov}
A. Akhmedov, J.B. Etnyre, T. E. Mark, I. Smith,
A note on Stein fillings of contact manifolds,
Math. Res. Lett. 15 (2008), no. 6, 1127--1132.

\bibitem{Auroux}
D. Auroux, V. Kulikov, V. Shevchishin,
Regular Homotopy of Hurwitz Curves,
Izv. Math. 68 (2004), 521--542.

\bibitem{AurouxSmith}
D. Auroux, I. Smith,
Lefschetz pencils, branched covers and symplectic invariants.
Symplectic 4-manifolds and algebraic surfaces,
Lecture Notes in Math. 1938, Springer, Berlin (2008), 1--53.

\bibitem{Bennequin}
D. Bennequin,
Entrelacements et \'equations de Pfaff,
Ast\'erisque 107-8 (1983), 87--161.

\bibitem{BirmanMenasco}
J. Birman, W. Menasco,
Stabilization in the braid groups II: Transversal simplicity of knots,
Geom. and Topol. 10 (2006), 1425--1452.

\bibitem{Ozbagci}
B. Ozbagci, A. I. Stipsicz,
Noncomplex smooth $4$-manifolds with genus-$2$ Lefschetz fibrations.
Proc. Amer. Math. Soc. 128 (2000), no. 10, 3125--3128.

\bibitem{Smith}
I. Smith,
Torus fibrations on symplectic four-manifolds.
Turkish J. Math. 25 (2001), no. 1, 69--95.

\bibitem{Teicher}
M. Teicher, M. Friedman,
On non fundamental group equivalent surfaces,
Alg. \& Geom. Topol. 8 (2008), 397--433.

\end{thebibliography}
\end{document}